\documentclass[twoside,10.5pt]{article}
\usepackage{mathrsfs}
\usepackage{pifont}
\usepackage{amsmath}
\usepackage{amsthm}
\usepackage{txfonts}
\usepackage{geometry}
\usepackage{latexsym}
\usepackage{amssymb}
\usepackage{graphicx}
\usepackage{geometry}
\usepackage{xcolor} 
\begin{document}
\newcommand{\arccot}{\mathop{\mathrm{arccot}}}
\newcommand{\artanh}{\mathop{\mathrm{artanh}}}
\newcommand{\arcoth}{\mathop{\mathrm{arcoth}}}
\newcommand{\li}{\mathop{\mathrm{li}}}
\newcommand{\const}{\mathop{\mathrm{const}}}

\begin{center}
{Originally appeared in: \\ 
{\it Prikladnaya Mekhanika i Matematika.} Sbornik Nauchnykh Trudov MFTI. Moscow, 1992, 96-103 (in Russian).}
\end{center}

\begin{center}
{\LARGE{Lyapunov Exponents for Burgers' Equation}}\\[20pt]
\end{center}

\begin{center}
Alexei Kourbatov$^1$

$^1$ Formerly with Moscow Institute of Physics and Technology \par
\smallskip
%
http://www.JavaScripter.net/math/pde/burgersequationdirichletproblem.htm  \\[25pt]
\end{center}



\textbf{Abstract}

We establish the existence, uniqueness, and stability of the stationary solution of
the one-dimensional viscous Burgers equation with the Dirichlet boundary conditions on a finite interval.
We obtain explicit formulas for solutions and analytically determine the Lyapunov exponents 
characterizing the asymptotic behavior of arbitrary solutions approaching the stationary one.



\medskip


\section*{Introduction}

Burgers' equation has the same nonlinearity form as the Navier-Stokes equations [1].
It is often used as a model equation in studying computational methods for solving 
partial differential equations (PDEs) [2].
In this paper we establish the existence, uniqueness, and stability of the stationary solution of
the  one-dimensional viscous Burgers equation (1) on a finite interval with the Dirichlet boundary conditions (4).
We use the Cole-Hopf transformation to give the result for any combination of $A$ and $B$ in the boundary conditions (4).
Using a different method (linearization) H.-O.~Kreiss and G.~Kreiss (1985) gave a similar result 
for a subset of cases: $A\ge|B|$, \ $B\le0<A$, as well as for Burgers' equation with forcing [8].
We obtain explicit formulas for solutions and analytically determine the Lyapunov exponents 
characterizing the asymptotic behavior of arbitrary solutions approaching the stationary solution
with the same boundary conditions (4).

\section{Explicit formulas for stationary solutions}

The viscous Burgers equation is the nonlinear partial differential equation
$$
u_t + u u_x = \nu u_{xx}  \eqno{(1)}
$$
with $\nu = \const > 0$.
If we set $u_t$ to zero, for the stationary solution $u=u^S(x)$ we obtain 
$$
u u_x = \nu u_{xx}.       \eqno{(2)}
$$
We note that $u u_x = \displaystyle{1\over2}(u^2)'_x$, therefore (2) gives
$$
2\nu u_x = u^2 + C_0.     \eqno{(3)}
$$

First, assume that $u_x \ne 0$ and $C_0$ is negative, $C_0=-a^2<0$ \ (i.\,e. $2\nu u_x < u^2$).
We have  $dx=2\nu du/(u^2-a^2)$,
$$
{ax\over\nu} = \ln \left( C_1 \left|{a-u\over a+u}\right|\right), \qquad\mbox{ where } \ 
C_1 = \left| {a+u(0)\over a-u(0)} \right|.
$$
If, in addition, $|u|<a$ \ (i.\,e. $u_x < 0$), then
$$
u = a\,{C_1-\exp(ax/\nu)\over C_1+\exp(ax/\nu)} = -2\nu k_0 \tanh(k_0(x-x_0)), \quad\mbox{ where } \ 
k_0 = {a\over2\nu}, \quad x_0 = 
                                {1\over k_0}\artanh{u(0)\over2\nu k_0},
$$
while if $|u|>a$ \ (i.\,e. $0 < 2\nu u_x < u^2$), then
$$
u = a\,{C_1+\exp(ax/\nu)\over C_1-\exp(ax/\nu)} = -2\nu k_0 \coth(k_0(x-x_0)), \quad\mbox{ where } \ 
k_0 = {a\over2\nu}, \quad x_0 = 
                                {1\over k_0}\arcoth{u(0)\over2\nu k_0}.
$$

Now assume that $C_0$ is positive, $C_0=a^2>0$ \ (i.\,e. $2\nu u_x > u^2$). Then $dx=2\nu du/(u^2+a^2)$,
$$
{ax\over2\nu} = \arctan {u\over a} +C_1,  \quad\mbox{ where } \ 
C_1 = -\arctan{u(0)\over a}; \qquad\mbox{ hence }\quad 
u = a \tan \left( {ax\over2\nu} - C_1 \right); 
$$
or, equivalently,
$$
u = -2\nu k_0 \cot(k_0(x-x_0)), \qquad\mbox{ where } 
k_0 = {a\over2\nu}, \quad x_0 = {1\over k_0} \arccot {u(0)\over2\nu k_0}.
$$

Finally, if \ $C_0 = 0$, then $u=-2\nu/(x-x_0)$; \ if $u_x=0$, then $u=\const$ and $u^2=|C_0|=\const$.
For convenience, all explicit formulas for stationary solutions are listed together in Table 1 (left column).

\vspace{4mm}
Table 1. Stationary solutions $u^S$ of Burgers equation and 
the corresponding solutions $\varphi^S$ of the heat equation (6). \linebreak
\phantom{Table 1.} 
$H=2\nu(B-A)-lAB$; \ \ 
$2\nu k_0=|C_0|^{1/2}$, where $C_0$ is the constant in (3); \ \ 
$u^S(x) = -2\nu(\ln|\varphi^S(x,t)|)'_x$. \ \ 

\vspace{-3mm}
\begin{center}
\begin{tabular}{lccl}\hline \phantom{\large$11^1$}
 Solution $u^S(x)$ of (1)
 &Conditions on  $u$, $u_x$
 &Conditions on $A$, $B$, $H$ 
 &Solution $\varphi^S(x,t)$ of (6)  \\[0.2ex]\hline
(a) \ $-2\nu k_0 \cot(k_0(x-x_0))$   &$0\le u^2 < 2\nu u_x$ &$A<B, \ H>0$ &$C\sin(k_0(x-x_0))\exp(-\nu k_0^2 t)$
                                                                                         \vphantom{\large$1^1$}  \\
\phantom{111.} \ $2\nu k_0 \tan(k_0(x-x^*_0))$ &\multicolumn{2}{c}{$\mspace{-9mu}$\it{(the same conditions and same solution as above)}}
                                                                          &$C\cos(k_0(x-x^*_0))\exp(-\nu k_0^2 t)$ \\
(b) \ $-2\nu/(x-x_0)$                &$0 <  u^2 = 2\nu u_x$ &$A<B, \ H=0$ &$C(x-x_0)$                            \\
(c) \ $-2\nu k_0\coth(k_0(x-x_0))$   &$0 < 2\nu u_x < u^2 $ &$A<B, \ H<0$ &$C\sinh(k_0(x-x_0))\exp(\nu k_0^2 t)$ \\
(d) \ $\pm2\nu k_0 = \const $        &$u_x = 0$             &$A=B$        &$C\exp(\nu k_0^2 t \mp k_0 x)$        \\
(e) \ $-2\nu k_0 \tanh(k_0(x-x_0))$  &$u_x < 0$             &$A>B$        &$C\cosh(k_0(x-x_0))\exp(\nu k_0^2 t)$ \\[0.2ex]\hline
\end{tabular}
\end{center} 
\vspace{2mm}

We will now consider Burgers equation (1) with the Dirichlet boundary conditions on the interval $x\in[0,l]$:
$$
u(0,t)=A, \qquad u(l,t)=B,   \eqno{(4)}
$$
where $A$ and $B$ are constants. 
Let us find out which explicit formulas (Table 1) can represent the stationary solution $u^S$ 
of equation (1) with boundary conditions (4). Here we are concerned exclusively with solutions that are 
continuous, bounded, and sufficiently smooth everywhere on the interval $x\in[0,l]$.

Clearly, when $A>B$, the stationary solution $u^S$ can only have form (e)
$-2\nu k_0 \tanh k_0(x-x_0)$ which is the only decreasing function in the left column of Table 1.
When $A=B$, the stationary solution $u^S$ can only have form (d);
all other explicit formulas for $u^S$ defined on $[0,l]$ are either strictly decreasing or 
strictly increasing functions of $x$.

To examine the stationary solution $u^S(x)$ for the trickiest case, $A<B$, we introduce the quantity
$$
H = 2\nu(B-A)-lAB.
$$
An elementary calculation shows that $u^S$ has form (b) if and only if $H=0$, $A<B$.
It remains to analyze the situations that yield solutions (a) and (c).
We note that, at any given point $(x,u(x))$, any graph $u^S$ of form (a) is steeper than (b),
while any graph of form (c) is less steep than (b).
Indeed, for any stationary solution $u^S$ we have a constant value of $C_0 = 2\nu u_x - u^2$;
solutions (a) are obtained from (3) when $2\nu u_x > u^2$ (steeper graphs, $C_0>0$, $H>0$),
while solutions (c) are obtained from (3) when $2\nu u_x < u^2$ (less steep graphs, $C_0<0$, $H<0$).
Thus when $A<B$ and $H>0$, we can only have $u^S$ given by formula (a);
when $A<B$ and $H<0$ we can only have $u^S$ given by formula (c).

Note also that we have not yet proved that a stationary solution satisfying boundary conditions (4)
{\it exists} for an arbitrary combination of $A$ and $B$. (We will prove this in Section 3.)
Still, in the simple cases (b) $A<B$, $H=0$ and (d)~$A=B$, it is already obvious that such stationary solutions do exist.

\section{The Cole-Hopf transformation}

Burgers equation (1) is a rare example of a nonlinear PDE that can be linearized 
using a simple transformation. Specifically, if in equation (1) we substitute
$$
u(x,t) = -2\nu(\ln|\varphi(x,t)|)'_x  \eqno{(5)}
$$
then for the unknown function $\varphi(x,t)$ we obtain the {\em heat equation}
$$
\varphi_t = \nu\varphi_{xx}.          \eqno{(6)}
$$
The substitution (5) is known as the Cole-Hopf transformation [1, 2, 5, 6].
Let us discuss some interesting properties of this transformation.

Firstly, transformation (5) can produce the same solution $u(x,t)$ of (1) from many different solutions $\varphi(x,t)$ of (6);
these $\varphi(x,t)$ may differ from each other by an arbitrary nonzero mutiplier $C$.
Indeed, $(\ln|\varphi|)'_x=(\ln|C\varphi|)'_x$ for any constant $C\ne0$.

Secondly, zero values of $\varphi(x,t)$ are mapped by (5) into discontinuities of $u(x,t)$.
Therefore, to get a continuous $u(x,t)$, it is not enough to start from a continuous solution 
$\varphi(x,t)$ of (6). We, moreover, need to restrict ourselves to those solutions $\varphi(x,t)$ 
that are nonzero everywhere on $[0,l]$ for all $t\ge0$.

Further, stationary solutions $u^S(x)$ of (1) correspond to solutions $\varphi^S(x,t)$ of (6)
that may or may not be stationary. 
Explicit formulas for those $\varphi^S(x,t)$ that yield stationary solutions $u^S(x)$ 
are listed in the right column of Table 1.
Interestingly, among these $\varphi^S(x,t)$ we find ``non-physical'' solutions of the heat equation
that grow infinitely large when $t\to\infty$.

\section{Existence and uniqueness of the stationary solution}

Using the Cole-Hopf transformation (5),
we will now establish the existence and uniqueness of the stationary solution of (1), (4) for any $A$ and $B$.
Note that (5) transforms the problem (1), (4) into the following problem for heat equation (6) 
with the Robin boundary conditions:
$$
\varphi_t = \nu\varphi_{xx}
$$
$$
\varphi_x(0,t) + {A\over2\nu}\varphi(0,t) = 0, \qquad
\varphi_x(l,t) + {B\over2\nu}\varphi(l,t) = 0. \eqno{(7)}
$$
Denote by $\varphi^S$ the solution of (6) that under transformation (5) yields the stationary solution $u^S$ of (1). 
Our $\varphi^S$ must have the form $\varphi^S(x,t)=X(x)\cdot T(t)$.
(This can be checked directly by substituting $\varphi^S$ into (5), or simply by inspection of the right column in Table 1.)
Here $X(x)$ is a function of the $x$ coordinate only, and $T(t)$ is a function of time $t$ only.
Substituting this $\varphi^S$ into the heat equation(6) and dividing through by $\nu T X$, we get
$$
{T'\over\nu T} = {X''\over X} = -\lambda.
$$
(One ratio is a function of $t$ only, while the other ratio is a function of $x$ only.
In order for these two ratios to be equal, they both must be equal to a constant which we denote $-\lambda$.)

For the function $X(x)$, problem (6), (7) translates into an eigenvalue problem (a Sturm-Liouville problem) 
with Robin boundary conditions:
$$
-X''(x) = \lambda X(x)       \eqno{(8)}
$$
$$
X'(0) + {A\over2\nu}X(0) = 0, \qquad
X'(l) + {B\over2\nu}X(l) = 0; \eqno{(9)}
$$
and for the function $T(t)$ we readily obtain 
$$
T(t) = C\exp(-\nu\lambda t).  \eqno{(10)}
$$

For $u^S$ to be continuous, $\varphi^S$ must be nonzero everywhere on the interval $[0,l]$.
So the question now is: how many eigenfunctions of (8), (9) are nonzero everywhere on $[0,l]$?
The answer is well known: for any $A$ and $B$, there is one and only one such eigenfunction.
This follows from the familiar fact that, for any $A$ and $B$ in problem (8), (9), all eigenvalues $\lambda_i$
$(\lambda_0 < \lambda_1 < \ldots)$ have multiplicity 1, and the respective eigenfunction $X_i(x)$
has exactly $i$ zeros inside the interval $(0,l)$; see [3, pp.\,14-18].
Thus, in problem (8), (9) we are interested in the eigenfunction $X_0(x)$ that has no zeros for $x\in[0,l]$ 
and corresponds to the least eigenvalue $\lambda_0$. For $\varphi^S$ we find, up to a nonzero multiplier $C$,
$$
\varphi^S(x,t) = C X_0(x)\cdot\exp(-\nu\lambda_0 t)\qquad\mbox{($\varphi^S$ has no zeros for $x\in[0,l]$).} 
$$

Therefore, for any $A$ and $B$, {\it there exists a unique stationary solution} $u^S(x)$ of Burgers equation (1)
with boundary conditions (4):
$$
u^S(x) = -2\nu(\ln|\varphi^S(x,t)|)'_x  = -2\nu(\ln|X_0(x)|)'_x.
$$

\section{Stability and Lyapunov exponents}

Now let us study the evolution of the absolute value $|u-u^S|$ for an arbitrary non-stationary solution
$$
u(x,t) = -2\nu(\ln|\varphi(x,t)|)'_x = -2\nu{\varphi_x(x,t)\over\varphi(x,t)},
$$
where both $u(x,t)$ and $u^S(x)$ satisfy the Burgers equation (1) with boundary conditions (4), and 
$\varphi(x,t)$ is a suitable positive solution of (6). It is known that
the solution $u(x,t)$ exists for ``reasonable'' combinations of the boundary conditions (4) and initial condition $u(x,0)$ [9].
We say that $u^S$ is {\it stable} if $|u-u^S|\to0$ as $t\to\infty$, for an arbitrary $u$ obeying (1), (4). We have
\begin{eqnarray*}
|u-u^S|
 &=&2\nu \left|{\varphi^S_x\over\varphi^S} - {\varphi_x\over\varphi} \right|
 ~=~ 2\nu \left|{\varphi\varphi^S_x - \varphi^S\varphi_x\over\varphi^S\varphi} \right|
\\
 &=&2\nu \left|{\varphi(\varphi^S_x-\varphi_x) + \varphi_x(\varphi-\varphi^S)\over\varphi^S\varphi} \right|
 ~=~ 2\nu \left|{\varphi_x\tilde{\varphi} - \varphi\tilde{\varphi}_x\over\varphi^S\varphi} \right|.
\end{eqnarray*}
Here we have introduced the notation $\tilde{\varphi}=\varphi-\varphi^S$.
Taking into account that $u=-2\nu\varphi_x/\varphi$, for all $x\in[0,l]$ and all $t\ge0$ we obtain the estimate
$$
|u-u^S| ~\le~ |u|\cdot\left|{\tilde\varphi\over\varphi^S} \right|
             + 2\nu \left|{\tilde\varphi_x\over\varphi^S} \right|
        ~\le~ \max_{x\,\in\,[0,l]}|u(x,0)| \cdot \left|{\tilde\varphi\over\varphi^S} \right|
             + 2\nu \left|{\tilde\varphi_x\over\varphi^S} \right|.
\eqno{(11)}
$$
In inequality (11) we have used the maximum principle for Burgers equation:
the solution $u(x,t)$ attains its maximum either in the initial value $u(x,0)$
or at the boundary of the interval $[0,l]$. 
(A discussion of maximum principles for PDEs can be found in [4, 7, 9].
The proof of the maximum principle for Burgers equation is similar to that for linear parabolic PDEs.)

Expand $\varphi(x,t)$ in a series over the system of eigenfunctions $X_i(x)$ of (8), (9):
$$
\varphi(x,t) = \sum_{i=0}^\infty \alpha_i X_i(x)T_i(t)
             = \sum_{i=0}^\infty \varphi_i(x,t),
\quad T_i(t) = \exp(-\nu\lambda_i t),
\quad T_i(0) = 1.
\eqno{(12)}
$$
In this series, the term $\varphi_0(x,t) = \alpha_0 X_0(x)T_0(t)$
is the same as $\varphi^S$ (Table 1) up to a constant nonzero multiplier.
Let us choose $C$ in the expression of $\varphi^S$ (Table 1) so that $\varphi_0 = \varphi^S$.
If we now compute the difference $\varphi-\varphi^S$, the term $\varphi_0(x,t)$ will cancel out, and we get
$$
\tilde\varphi \,=\, \varphi-\varphi^S = \sum_{i=1}^\infty \varphi_i(x,t).
\eqno{(13)}
$$
Since $T_i(t) = \exp(-\nu\lambda_i t)$, we see that $\varphi_1(x,t)$ becomes the {\it largest} term in (13) when $t\to\infty$ 
(assuming $\alpha_1\ne0$ in (12)). We then have 
$$
\max_{x\,\in\,[0,l]}|\varphi^S      | \asymp \exp(-\nu\lambda_0 t), \quad
\max_{x\,\in\,[0,l]}|\tilde\varphi  | \asymp \exp(-\nu\lambda_1 t), \quad
\max_{x\,\in\,[0,l]}|\tilde\varphi_x| \asymp \exp(-\nu\lambda_1 t)  \qquad\mbox{ as } t\to\infty,
$$
so the estimate (11) results in 
$$
\max_{x\,\in\,[0,l]} |u - u^S| \,\asymp\, \exp(-\nu(\lambda_1-\lambda_0) t) \qquad\mbox{ as } t\to\infty.
\eqno{(14)}
$$
This paves the way to proving the {\it stability of the stationary solution} $u^S$. 
Indeed, the difference $|u - u^S|$ is an exponentially vanishing quantity when $t\to\infty$.
Nevertheless, the convergence of $|u - u^S|$ to zero might turn out to be very slow; this is the case when
the least two eigenvalues $\lambda_0$ and $\lambda_1$ in problem (8), (9) differ only slightly. 

We got the estimate (14) under the assumption that $\alpha_1\ne0$ in (12), that is,
in the series expansion of $\varphi$ over the system of eigenfunctions $X_i(x)$ 
there is a nonzero term $\varphi_1$ containing the eigenfunction $X_1(x)$.
However, if it so happens that one or more initial terms in (13) are zero,
then the series (13) for $\tilde\varphi = \varphi-\varphi^S$ will start at
some $\varphi_n$ ($n>1$). In the general case, therefore, instead of (14) we would have
$$
\max_{x\,\in\,[0,l]} |u - u^S| \,\asymp\, \exp(-\nu(\lambda_n-\lambda_0) t) \qquad\mbox{ as } t\to\infty,
\eqno{(15)}
$$
where $n$ is the number of the first nonzero term in the series expansion of $\tilde\varphi=\varphi-\varphi^S$ (13).
We have thus proved that the stationary solution $u^S$ is stable:
$|u - u^S|\to0$ as $t\to\infty$.

Note that the functions $\varphi_i$ $(i=1,2,\ldots)$ in (12) 
have the same explicit formulas as $\varphi^S$ (Table 1),
except that each $\varphi_i$ contains its own values in place of $k_0$ and $x_0$; 
let us denote these new constant values by $k_i$ and $x_i$, respectively.

All constants $k_i$ and $x_i$ can be found if we substitute the general solutions of (8) 
(trigonometric, exponential or hyperbolic functions) 
for the eigenfunctions $X_i(x)$ $(i=0,1,2,\ldots)$ in the boundary conditions (9). 
In most cases (i.\,e., cases (a), (c), (e) in Table 1), this substitution yields the following transcendent equations 
for $\xi_i = k_i l$:
$$
\cot\xi_i ~=~ {p\over\xi_i} ~+~ q\xi_i \quad\mbox{ for } \varphi_i \mbox{ of form (a) in Table 1, }
  \quad\lambda_i=k_i^2>0, 
  \quad\mbox{ and }   \eqno{(16)}
$$
$$
\coth\xi_i ~=~ {p\over\xi_i} ~-~ q\xi_i \quad\mbox{ for } \varphi_i \mbox{ of form (c) or (e) in Table 1, }
  \quad \ \ \lambda_i=-k_i^2<0,         
  \eqno{(17)}
$$
$$
\mbox{ where }
\quad \xi_i = k_i l > 0,
\quad p = {lAB \over 2\nu(B-A)},
\quad q = {2\nu \over l(B-A)}.
$$
The transcendent equation (16), with $\cot\xi_i$, may correspond to {\it any} $i$, whereas
equation (17), with $\coth\xi_i$, may correspond only to $i=0, 1$ 
(the least two eigenvalues $\lambda_0$, $\lambda_1$)
because hyperbolic functions cannot have more than one zero value on the interval $[0,l]$.

When $A=B$ in (4) and (9), we have an exceptional case: all $k_i$ and $\lambda_i$ can be found in a closed form.
Here the interval $[0,l]$ contains a whole number of semiperiods of the eigenfunction 
$X_i(x) = \sin(k_i(x-x_i))$, $i=1,2,\ldots$\,, which readily yields
$$
k_i = {\pi i \over l} \ \ (i=1,2,\ldots), \quad\mbox{ while }\quad k_0={|A|\over2\nu},
\quad X_0(x) = C\exp(\pm k_0 x); \ \ 
\mbox{ see Table 1 (d).}
$$ 
Therefore, if $A=B$, we find
$$
      \lambda_n = \left( {\pi n \over l} \right)^2 \ \ (n\ge1),
\quad \lambda_0 = -\left( {A\over2\nu} \right)^2,
\quad\mbox{ and }
\quad \lambda_n-\lambda_0 ~=~ \left( {\pi n \over l} \right)^2 +\, \left( {A\over2\nu} \right)^2;
\quad\mbox{cf.\ (14), (15).}
$$

Now we will reuse the customary definition of {\it Lyapunov exponents} 
in the context of problem (1), (4) for Burgers equation. 
Let $u(x,t)$ be a solution of (1),(4). The Lyapunov exponent $\mu$ of this solution is defined as
$$
\mu = \limsup_{t\to\infty} {\ln|| u - u^S || \over t}.
\eqno{(18)}
$$
This definition, in general, depends on our choice of the norm $||\cdot||$.
If $u(x,t)$ behaves so that  $|| u - u^S || \asymp \exp(\delta t)$ as $t\to\infty$, 
then it is easy to see that $\delta$ is the Lyapunov exponent of this $u(x,t)$.

Let us use the norm defined as the maximum absolute value:
$$
||w(x)|| = \max_{x\,\in\,[0,l]} |w(x)|.
$$
Then estimates (14), (15) allow us to determine all Lyapunov exponents for any $u(x,t)$ satisfying (1), (4):
$$
\mu_i = -\nu(\lambda_i-\lambda_0), \quad i=1,2,\ldots, 
\eqno{(19)}
$$
where, as before, $\lambda_i$ are eigenvalues of (8), (9). 
Solutions $u(x,t)$ corresponding to the Lyapunov exponents $\mu_i$ can be written simply as
$$
u_i(x,t) = -2\nu(\ln|\varphi^S(x,t)+\varphi_i(x,t)|)'_x, \quad i=1,2,\ldots,
$$
where $\varphi_i(x,t)$ is the respective term of (12). For example, when $u^S$ has the form (a) in Table 1,
we have
\begin{eqnarray*}
\varphi^S(x,t) &=& C\sin(k_0(x-x_0))\exp(-\nu k_0^2 t) 
                   \qquad\mbox{($\varphi^S$ has no zeros for $x\in[0,l]$),}  \\
\varphi_i(x,t) &=& \alpha_i\sin(k_i(x-x_i))\exp(-\nu k_i^2 t) 
                   \qquad\mbox{ ($\varphi_i$ has $i$ zeros for $x\in[0,l]$), }  
\end{eqnarray*}
and we can write a solution $u_i(x,t)$ corresponding to the Lyapunov exponent $\mu_i$ as follows:
$$
u_i(x,t) = -2\nu {
Ck_0\cos(k_0(x-x_0)) + \alpha_i k_i\cos(k_i(x-x_i))\cdot\exp(-\nu(k_i^2 - k_0^2) t)
\over
C\sin(k_0(x-x_0)) + \alpha_i \sin(k_i(x-x_i))\cdot\exp(-\nu(k_i^2 - k_0^2) t)
}.
\eqno{(20)}
$$
Because each individual term in series (12) satisfies the Robin boundary conditions (7),
each function $u_i(x,t)$ defined as above must satisfy the Dirichlet boundary conditions (4).

We have thus determined the Lyapunov exponents in the nonlinear problem (1), (4) for 
Burgers equation: we have found that formula (19) relates the Lyapunov exponents $\mu_i$ 
to the eigenvalues $\lambda_i$ of the linear problem (8), (9). 
All Lyapunov exponents $\mu_i$ are negative; there are countably many of them; we can write 
explicit formulas for the corresponding solutions $u_i(x,t)$ of Burgers equation (1).
This is an interesting example of a situation where one can analytically determine 
the Lyapunov exponents for solutions of a nonlinear PDE with Dirichlet boundary conditions. 

\bigskip

\smallskip\noindent
\textbf{References} 

\hangafter=1
\setlength{\hangindent}{2em}
[1] Karpman, V.~I.~(1975). {\it Nonlinear Waves in Dispersive Media}. International Series of Monographs in Natural Philosophy, Vol.\,71,
Pergamon, 1975.

\hangafter=1
\setlength{\hangindent}{2em}
[2] Fletcher, C.~A.~J.~(1991). {\it Computational Techniques for Fluid Dynamics}. Springer, 1991.

\hangafter=1
\setlength{\hangindent}{2em}
[3] Levitan, B.~M.,~and Sargsjan, I.~S.~(1975). {\it Introduction to Spectral Theory}. 
Translations of Mathematical Monographs, Vol.\,39, AMS, 1975.

\hangafter=1
\setlength{\hangindent}{2em}
[4] Vladimirov, V.~S.~(1984). {\it Equations of Mathematical Physics}. Mir Publishers, 1984.

\smallskip\noindent
{\it References added in the English version:}

\hangafter=1
\setlength{\hangindent}{2em}
[5] Hopf, E. (1950). The partial differential equation $u_t+uu_x=\mu u_{xx}$. 
{\it Comm.~Pure and Appl.~Math.} {\it 3}, 201-230.

\hangafter=1
\setlength{\hangindent}{2em}
[6] Cole, J.~D.~(1951). On a quasilinear parabolic equation occurring in aerodynamics, 
{\it Quart. Appl. Maths.} {\it 9}, 225-236.

\hangafter=1
\setlength{\hangindent}{2em}
[7] Protter, M.~H.,~and Weinberger, H.~F.~(1984). {\it Maximum Principles in Differential Equations}. Springer, 1984.

\hangafter=1
\setlength{\hangindent}{2em}
[8] Kreiss, H.-O., and Kreiss, G.~(1985)
Convergence to steady state of solutions of Burgers' equation. {\it NASA Contractor Report 178017}, ICASE No.\,85-50.
NASA Langley Research Center: Institute for Computer Applications in Science and Engineering. Hampton, VA, 1985.

\hangafter=1
\setlength{\hangindent}{2em}
[9] Ladyzhenskaja, O.~A., Solonnikov, V.~A., and Uralceva, N.~N.~(1967). 
{\it Linear and Quasilinear Equations of Parabolic Type}, Nauka, Moscow, 1967 (in Russian).
English translation: AMS, Providence RI, 1968.


\vspace{0.5cm}
\textbf{Copyrights}

Copyright for this article is retained by the author.
\end{document}